\documentclass[12pt]{amsart}

\usepackage{graphicx}

\usepackage[english]{babel}
\usepackage{amsfonts}
\usepackage{amscd}
\usepackage{amsmath}

\newtheorem{Thm}{Theorem}
\newtheorem{rmk}{Remark}
\newenvironment{Rmk}{\begin{rmk}\em}{\end{rmk}}
\newtheorem{exm}{Example}

\newtheorem{prf}{Proof}

\newenvironment{Prf}{\begin{prf}\em}{\qed\end{prf}}
\newtheorem{prff}{}

\newtheorem{stat}{}

\newcommand{\NOT}[1]{}

\newcommand{\pa}{\par\medskip}

\newcommand{\LP}{\left(}  \newcommand{\RP}{\right)}

\newcommand{\BE}{\begin{equation}}  \newcommand{\EE}{\end{equation}}
\newcommand{\BR}{\begin{eqnarray*}}  \newcommand{\ER}{\end{eqnarray*}}
\newcommand{\BER}{$$\begin{array}}  \newcommand{\EER}{\end{array}$$}

\newcommand{\I}{\infty}

\newcommand{\al}{\alpha}
\newcommand{\be}{\beta}

\newcommand{\Om}{\Omega} 
 
\newcommand{\tfr}[2]{{\textstyle\frac#1#2}}

\newcommand{\cE}{\mathcal{E}}
\newcommand{\cF}{\mathcal{F}}

\newcommand{\cI}{\mathcal{I}}

\newcommand{\cO}{\mathcal{O}}

\newcommand{\cU}{\mathcal{U}}

\newcommand{\half}{\tfr12}

\title[No $O(N)$ queries]{No $O(N)$ queries for checking if $N$ intervals cover everything or for piercing $N$ pairs of intervals. An $O(N\log N)$-steps algorithm for piercing}

\author{Meir Katchalski}
\address{Department of Mathematics,
Technion -- Israel Institute of Technology,
Haifa 32000, Israel}
\email{meirk@techunix.technion.ac.il}
\author{Eliahu Levy}
\address{Department of Mathematics,
Technion -- Israel Institute of Technology,
Haifa 32000, Israel}
\email{eliahu@techunix.technion.ac.il}

\date{}

\begin{document}

\begin{abstract}
The complexity of two related geometrical (indeed, combinatorial) problems is considered, measured by the number of queries needed to determine the solution. It is proved that one cannot check in a linear in $N$ number of queries whether $N$ intervals cover a whole interval, or whether for $N$ pairs of intervals on two lines there is a pair of points intersecting each of these pairs of intervals (``piercing all pairs of intervals''). The proofs are related to examples which show that there is no ``Helly property'' here -- the whole set of $N$ may cover the whole interval (resp.\ may have no pair of points piercing all pairs of intervals) while any proper subset does not. Also, for the piercing problem we outline an algorithm, taking $O(N\log N)$ steps, to check whether there is a pair of points piercing all pairs of intervals and if there is, to find it.
\end{abstract}

\maketitle


\section{Introduction}
We consider two related settings. They are here formulated for rectangles (``two dimensions'') although they make sense for any number of dimensions. \pa

I. Consider $N$ rectangles $R_i=[a_i,b_i]\times[c_i,d_i],\,\,i=1\ldots,N$ which are subrectangles of a rectangle $R_0=[a_0,b_0]\times[c_0,d_0]$. Consider the sets of symbols:
\BR
&&X=\{a_0,a_1,\ldots,a_N,b_0,b_1,\ldots,b_N\},\\
&&Y=\{c_0,c_1,\ldots,c_N,d_0,d_1,\ldots,d_N\}.
\ER
One asks the question: do these subrectangles cover $R_0$? \pa

One is allowed queries of the form: for two elements $x,x'$ of $X$ (resp.\ $Y$), which of $x>x'$, $x<x'$, $x=x'$ does hold? Of course, we know in advance that:
$$a_0\le a_i<b_i\le b_0,\,i=1,\ldots,N,\quad
c_0\le c_i<d_i\le d_0,\,i=1,\ldots,N.\eqno(*)$$
And one wishes to determine, by some number of such queries, whether the $N$ rectangles cover $R_0$. In particular, can $O(N)$ queries suffice when $N\to\I$? \pa

II. Again $[a_i,b_i]$ and $[c_i,d_i]$ are subintervals of the sides of the rectangle $R_0=[a_0,b_0]\times[c_0,d_0]$. But now one considers ``crosses'' $C_i$ $:=$ the complement in $R_0$ of the Cartesian product of the complements in the sides of $R_0$ of $[a_i,b_i]$ and $[c_i,d_i]$, resp. And one wishes to determine, by queries as above, whether the {\em intersection} of the crosses is empty, i.e.\ whether the {\em union} of the Cartesian products of the complements is the whole $R_0$. \pa

This may be formulated as follows. View a point in $R_0$ as a pair of a point in $[a_0,b_0]$ and a point in $[c_0,d_0]$. We wish to check, by queries, whether the intersection of the ``crosses'' is {\em not} empty, thus whether there is a pair of points that belongs to all the crosses, i.e.\ meets either $[a_i,b_i]$ or $[c_i,d_i]$ for all $i=1,\ldots,N$. \pa

So our problem takes the form: on two fixed lines (now viewed as parallel), $N$ pairs of intervals, each pair consisting of one interval on each of two lines, are given. To determine by queries if there is a pair of points, one on each line, that have a common point with each pair of intervals (``can we pierce all pairs?''). Thus setting II will be referred to as the {\em piercing problem}. \pa

Now, by ``closing'' the sides of $R_0$ to circles (``at infinity''), thus making $R_0$ into a torus, one makes settings I and II look similar -- the complements of the ``crosses'', about whose union one asks, are then also a kind of ``rectangles'', but there is an important difference -- in setting II one asks about a union of ``rectangles'' which {\em all have a common point} (at ``infinity''). So setting II is, in fact, a restricted case of setting I. \pa

Note, that if by the answers to a set of queries we can determine the full ordering of $X$ and $Y$, then the whole ``geometrical picture'' is determined and hence so is the answer to the question about the union or intersection. But it is well-known that a set of $N$ elements can be ordered by $O(N\log N)$ queries (see the appendix). Therefore $O(N\log N)$ queries always suffice for our settings. \pa

But are $O(N)$ queries enough? We prove that:

\begin{Thm}\label{Thm:thm1}
For setting I, even in dimension $1$ (and consequently for higher dimensions), $O(N)$ queries do not suffice. Indeed, $\Om(N\log N)$ queries may be needed in the worst case.
\end{Thm}

For setting II in dimension $1$, one has to find whether the intersection of intervals is empty, which, as is well-known, is easily done by $O(N)$ queries (see Remark \ref{Rmk:intersection}).

\begin{Thm}\label{Thm:thm2}
For setting II in dimension $2$ (and consequently for higher dimensions), $O(N)$ queries do not suffice. Indeed, $\Om(N\log N)$ queries may be needed in the worst case.
\end{Thm}

For setting II (piercing), in dimension $2$, we shall also outline an algorithm taking $O(N\log N)$ steps that, provided we already know the ordering of $X$ and $Y$, determines whether the $N$ pairs of intervals can be pierced and, if they can, finds a piercing pair of points (i.e., determines whether the $N$ crosses intersect and, if they can, finds a point in the intersection). \pa

The proofs (and clarification of the notions) will be given in the following sections.

\section{Preliminaries to the proofs}
To fix notation, we here continue to speak about the two-dimensional case as in the introduction. \pa

Our knowledge after a sequence of queries may be described as a preorder $\al$ on $X$ and a preorder $\be$ on $Y$. \pa

Recall that a preorder is a reflexive and transitive relation $\le$. It defines an equivalence relation $x\equiv y:=x\le y\,\&\,y\le x$ and a partial ordering in the factor set. We call a preorder {\em full} if for any two elements $x$ and $y$ either $x\le y$ or $y\le x$, or both. $<$ will always mean $(\le\,\&\,\neg\ge)$. \pa

$\al$ will be the smallest preorder that contains the answers to the queries and the part of $(*)$ that refers to $X$, similarly $\be$ for $Y$. For any such preorder $\al$, denote by $\cO(\al)$ the set of all full preorders on $X$ that extend $\al$ (in other words, are strengthenings of $\al$, where in speaking of strengthening and weakening of a preorder we allow equality). similarly for $\cO(\be)$. Note that a full preorder $\xi$ on $X$ with a full preorder $\eta$ on $Y$ determine completely the geometric configuration of the rectangles. \pa

Denote by $\cF(X)$ (resp.\ $\cF(Y)$) the set of all full preorders on $X$ (resp.\ $Y$) and let $\cU$ be the subset of $\cF(X)\times\cF(Y)$ consisting of all pairs of full preorders in whose configuration of rectangles $\cup_{i=1}^NR_i=R_0$ (for setting I) or the intersection of the crosses is empty (for setting II). \pa

Now, if we can decide whether the union of rectangles is $R_0$ (or the intersection of the crosses is empty) by $m$ queries, then, since each query can have at most $6$ answers ($=,\ne,>,<,\ge,\le$), there are at most $6^m$ different possible sequences of answers to the $m$ queries, each such sequence holding precisely when our configuration $(\xi,\eta)$, $\xi\in\cF(X),\eta\in\cF(Y)$ is in some product set $\cO(\al)\times\cO(\be)$. $\cU$ must be the union of those products for which the queries give the answer ``yes'', so $\cU$ is a union of no more than $6^m$ products of the form $\cO(\al)\times\cO(\be)$. \pa

Put otherwise: the smallest number of queries that suffice is no less than the $\log_6$ of the smallest number $M$ of products $\cO(\al)\times\cO(\be)$ such that $\cU$ can be written as the union of them. \pa

\begin{Rmk} \label{Rmk:intersection}
For the problem whether the {\em intersection} of the rectangles is empty there is a well-known simple way to answer by $O(N)$ queries: ask whether $a_1>a_2$. If $a_1>a_2$ forget about $a_2$ and consider $a_1$. Now compare it with $a_3$ and forget the smaller, etc. Similarly with the $b$'s. After that check whether the remaining interval $[a,b]$ is empty. The same with the $c$'s and $d$'s. \pa

And indeed, the set $\cI\subset\cF(X)\times\cF(Y)$ of all the configurations with nonempty intersection of the rectangles is a single product set $\cO(\al)\times\cO(\be)$ with
\BR
\al&:&a_1,\ldots,a_N\le b_1,\ldots,b_N,\\
\be&:&c_1,\ldots,c_N\le d_1,\ldots,d_N.
\ER
\end{Rmk}

\section{Proof of Theorem \ref{Thm:thm1} (for dimension $1$)}
\begin{Prf}
Since we are considering dimension $1$, we will have only $X$ and only intervals (instead of rectangles) $R_i=[a_1,b_i]$ and $R_0=[a_0,b_0]$, and we wish to show that by less than $\Om(N\log N)$ queries we cannot always determine whether $\cup_{i=1}^NR_i=R_0$. \pa

Retain the notation of the previous section (adapted to setting I in dimension $1$). Thus $\cU$ is the subset of $\cF(X)$ consisting of all full preorders in whose configuration of intervals $\cup_{i=1}^NR_i=R_0$. \pa

Now suppose $\cU$ is the union of $M$ sets of the form $\cO(\al_j),\,j=1\ldots,M$ where $\al_j$ are preorders in $X$, and assume, as usual, that $\al_j$ contains (i.e.\ implies) $(*)$ for $X$. We shall show that $M$ must be something like $N!$. \pa

To this end, Let $P=(i_1,\ldots,i_N)$ be an ordering of $\{1,\ldots,N\}$, in other words, a permutation of $\{1,\ldots,N\}$. Let $\xi_P\in\cF(X), i=1\ldots N$ be the full preorder defined by:
\BE\begin{array}{ccc}
&&a_0=a_{i_1}<a_{i_2}<b_{i_1}<a_{i_3}<b_{i_2}<a_{i_4}<b_{i_3}<a_{i_5}<b_{i_4}<\ldots\\
&&\ldots<a_{i_N}<b_{i_{N-1}}<b_{i_N}=b_0,i=1,\ldots N.\\
\end{array}\label{eq:xi}
\EE
Here the intervals $R_i$ of $\xi_P$ form a kind of ``chain'': $R_1,\ldots,R_N$ are ordered as $P$, with each one partly overlapping the next. We have $\cup_{i=1}^NR_i=R_0$, hence $\xi_R\in\cU$. \pa

Therefore, there is an $\al_j$, which we shall denote by $\al_P$, such that $\xi_P\in\cO(\al_j)\subset\cU$. \pa

$\al_P$, as a preorder, is a weakening of $\xi_P$. I claim that $\al_P$ contains the inequality $a_{i_k}\le b_{i_{k-1}}$ for each $k=2,\ldots,N$. Because if not, than the full preorder $\xi'$ obtained from $\xi_P$ by flipping $a_{i_k}$ and $b_{i_{k-1}}$ would be a strengthening of $\al_P$, so $\xi'\in\cO(\al_P)$, contrary to the fact that $\xi'\notin\cU$ (two adjacent links in the chain are no longer overlapping, so the union fails to be $R_0$). \pa

Thus, for each link in the chain of strict inequalities in (\ref{eq:xi}), hence for any strict inequality deduced from them, $\al_P$ contains that inequality with $\le$ instead of $<$. \pa

Now let $P'=({i'}_1,\ldots,{i'}_N)$ be some other ordering of $\{1,\ldots,N\}$. I claim that $\al_P=\al_{P'}$ is possible for only if $P'=P$. \pa

Indeed, let $k\mapsto\ell=\ell(k)$ be the one-one mapping defined by $i_k={i'}_{\ell(k)}$. Suppose first that there exist $k,k'=1,\ldots,N$ such that $k'\ge k-1$ and $\ell(k')\le\ell(k)-2$. Then in $\xi_{P'}$,\,\, $a_{i_k}=a_{{i'}_{\ell(k)}}>b_{{i'}_{\ell(k)-2}}\ge b_{{i'}_{\ell(k')}}=b_{i_{k'}}$. Thus $\al_{P'}$, a weakening of $\xi(P')$, cannot contain $a_{i_k}\le b_{i_{k'}}$, which $\al_P$ contains (since $k'\ge k-1$), and one finds that in this case $\al_{P'}\ne\al_P$. \pa

This means that $\al_{P'}$ can equal $\al_P$ only if one has always
\BE\label{eq:kI}
k'\ge k-1\Rightarrow\ell(k')\ge\ell(k)-1.
\EE
This implies that for each $k$ the one-one mapping $k\mapsto\ell(k)$ must map the set $\{k'=1,\ldots,N\,|\,k'\ge k-1\}$ into $\{L'=1,\ldots,N\,|\,L'\ge\ell(k)-1\}$, which is possible only if $\ell(k)\le k$ for all $k=2,\ldots,N$, while $\ell(1)\le2$, and also, for each $k'$, must map the set $\{k=1,\ldots,N\,|\,k'\ge k-1\}$ into $\{L=1,\ldots,N\,|\,\ell(k')\ge L-1\}$, which is possible only if $k'\le\ell(k')$ for $k'=1,\ldots,N-1$ while $\ell(N)\ge N-1$. Assuming $N\ge3$ one deduces that $\ell(k)=k$ for all $k$, hence $k\mapsto\ell(k)$ is the identity mapping and $P=P'$. \pa

We have thus found that all the $\al_P$ for the $N!$\,\,$P$'s are different. Therefore $M\ge N!$ and the number of queries must be
$$m\ge\log_6(N!)=\Om(N\log N).$$
\end{Prf}
\begin{Rmk}\label{Rmk:rmk1}
The ``chain'' in the proof above shows trivially that there can be $N$ intervals whose union is $R_0$ while for any proper subset of them the union is not $R_0$ (compare Remark \ref{Rmk:rmk2}).
\end{Rmk}
\begin{Rmk}
This problem whether the UNION of intervals covers everything (setting I for dimension $1$), is related to the EQUALITY checking problem: suppose $\{x_1,\ldots,x_N\}$ are $N$ unknown numbers, we are allowed queries of the form: is $x_i$ greater, less or equal to $x_j$? and one wishes to determine by such queries whether all the $N$ numbers are different or there is some equality among them. It seems to be well-known that this cannot be done by less than $O(N\log N)$ queries (clearly $O(N\log N)$ queries suffice, since they suffice to determine the full order of the $x$'s).
\begin{itemize}
\item[1.]
One might say that this proves the claim of Theorem \ref{Thm:thm1}: if one could solve UNION in less than $O(N\log N)$ queries then one could solve EQUALITY also -- think of each $x_i$ as an interval of the form $[m_i,m_i+1]$, the $m_i$ being unknown integers in $\{0,1,\ldots,N-1\}$. The union covers $[0,N]$ iff these intervals are all distinct. \pa
 
A counter to this is that with UNION we are allowed different queries for the ends of the interval (here $m_i$ and $m_i+1$), which might make us perform better. \pa
\item[2.]
The claim that $O(N\log N)$ queries are needed for EQUALITY might seem to be almost obvious:
if one has a number of queries which is less then what is needed to order $N$ numbers by comparisons then there are two numbers $x_i$ and $x_j$ for which one does not know the relation between them so certainly one does not know if they are equal. It is very well known that ordering $N$ numbers requires $O(N\log N)$ queries, so EQUALITY needs $O(N\log N)$ too. \pa

Yet, this argument might seem too heuristic. After all, that ordering needs $O(N\log N)$ queries follows from a consideration of the information needed: there are $N!$ orderings, $m$ queries give
$\exp(O(m))$ possible answers, so we need $\log(N!)=O(N\log N)$ queries. For EQUALITY the answer is just one bit. \pa

\item[3.]
Anyhow, making the argument of 2. more precise leads us to a proof that EQUALITY needs $O(N\log N)$ queries which follows the lines of the above proof of Theorem \ref{Thm:thm1} for UNION (in a simpler way). \pa

Indeed, one has to prove that if one writes the set $\cU$ of all full preorders of $\{1,2,...,N\}$ that are permutations (i.e.\ no equality) as a union of $\cO(\al_j)$, $j=1,...,M$, $\al_j$ preorders, then $M$ is at least something like $N!$. \pa

Now let $P=(i_1,...,i_N)$ be a permutation of $\{1,...,N\}$, and consider $P$ as a full preorder. Since it is in $\cU$, it is in some $\al_i$, which we denote by $\al_P$. Since $\al_P$ is a strengthening of $P$, it must allow $i_{k-1}<i_k$. If it allows more (i.e.\ $i_{k-1}<=i_k)$ then define a full preorder $\xi'$ by changing $i_{k-1}<i_k$ to $i_{k-1}=i_k$ in $P$. Then $\xi'$ is not in $\cU$ but is a strengthening of $\al_P$, thus is in $\cO(\al_P)$ -- impossible. Therefore $\al_P$ requires $i_{k-1}<i_k$, and that for all $k=2,...,N$, making $\al_P=P$. Consequently, $P\mapsto\al_P$ is one-one, hence $M\ge N!$.
\end{itemize}
\end{Rmk}

\section{Proof of Theorem \ref{Thm:thm2} (for dimension $2$)}
\begin{Prf}
The proof proceeds similarly to the proof of Theorem \ref{Thm:thm1}. The main difference is in the configuration (now two-dimensional) that one corresponds to a permutation $P$.
\pa

We retain the notation of the previous sections (adapted to setting II in dimension $2$). Thus $\cU$ is the subset of $\cF(X)\times\cF(Y)$ consisting of all pairs of full preorders in $X$ and $Y$, resp.\ in whose configuration of intervals the intersection of the crosses is empty. \pa

Now we suppose $\cU$ is the union of $M$ sets of the form $\cO(\al)\times\cO(\be)$
where $\al_j$ and $\be_j$ are preorders in $X$ and $Y$, resp.\ and as usual we assume that $\al_j$ and $\be_j$ contain (i.e.\ imply) $(*)$ for $X$ and $Y$. We shall show that $M$ must be at least something like $\tfr14\LP\lfloor\half N\rfloor!\RP^2$. \pa

To this end, let $P=(i_1,\ldots,i_N)$ be an ordering of $\{1,\ldots,N\}$, in other words, a permutation of $\{1,\ldots,N\}$. But here we assume that $P$ preserves the parity of indices, i.e.\ that $i_k$ is odd for odd $k$ and even for even $k$. To $P$ we correspond full preorders $\xi_P\in\cF(X)$ and $\eta_P\in\cF(Y)$, that will make the ``crosses'' (here L-shaped) look as in Figure \ref{fig:fig1} for $N$ even or Figure \ref{fig:fig2} for $N$ odd. \pa

Here, for each ``upper'' or ``lower'' right angle of the ``stairs'' there is an L-shaped ``cross'' with the inner sides of arms of the ``L'' extending the sides of that angle. The ordering of the crosses along the stairs (alternatingly with inner angle ``upward'' and ``downward'') will follow $P$. Thus the intervals of the cross extending the lowest right angle of the stairs are $[a_{i_1},b_{i_1}]$ and $[c_{i_1},d_{i_1}]$, those of the next cross, extending the next right angle of the stairs, are $[a_{i_2},b_{i_2}]$ and $[c_{i_2},d_{i_2}]$ etc. \pa
\newpage
Thus, the full preorders $\xi_P$ and $\eta_P$ for $N=8$ are given by (recall that the two intervals which define a ``cross'' span the ``hub'' of the cross)
\begin{figure}[h]
\includegraphics[width=0.24\textwidth]{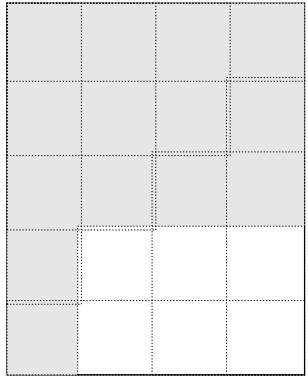}
\caption{``Staircase'' for $N=8$,
the third ``cross'' shaded.}
\label{fig:fig1}\end{figure}
\BR
&&a_0=a_{i_1}=a_{i_3}=a_{i_5}=a_{i_7}=\\
&&=b_{i_1}<b_{i_3}<a_{i_2}<b_{i_5}<a_{i_4}<b_{i_7}<a_{i_6}<a_{i_8}=\\
&&=b_{i_2}=b_{i_4}=b_{i_6}=b_{i_8}=b_0,\\
&&c_0=c_{i_2}=c_{i_4}=c_{i_6}=c_{i_8}=\\
&&=d_{i_2}<c_{i_1}<d_{i_4}<c_{i_3}<d_{i_6}<c_{i_5}<d_{i_8}<c_{i_7}<\\
&&<d_{i_1}=d_{i_3}=d_{i_5}=d_{i_7}=d_0.
\ER
And for $N=9$\,\,$\xi_P$ and $\eta_P$ are given by
\begin{figure}[h]
\includegraphics[width=0.3\textwidth]{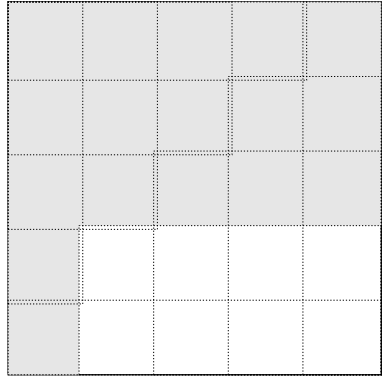}
\caption{``Staircase'' for $N=9$,
the third ``cross'' shaded.}
\label{fig:fig2}\end{figure}
\BR
&&a_0=a_{i_1}=a_{i_3}=a_{i_5}=a_{i_7}=a_{i_9}=\\
&&=b_{i_1}<b_{i_3}<a_{i_2}<b_{i_5}<a_{i_4}<b_{i_7}<a_{i_6}<b_{i_9}<a_{i_8}<\\
&&<b_{i_2}=b_{i_4}=b_{i_6}=b_{i_8}=b_0,\\
&&c_0=c_{i_2}=c_{i_4}=c_{i_6}=c_{i_8}=\\
&&=d_{i_2}<c_{i_1}<d_{i_4}<c_{i_3}<d_{i_6}<c_{i_5}<d_{i_8}<c_{i_7}<c_{i_9}=\\
&&=d_{i_1}=d_{i_3}=d_{i_5}=d_{i_7}=d_{i_9}=d_0.
\ER

The rest of the proof largely imitates the corresponding part of the proof of Theorem \ref{Thm:thm1}. \pa

The intersection of the crosses of $\xi_P$ and $\eta_P$ is empty (otherwise put: the union of the rectangles that are the complements of the crosses is the whole $R_0$). Therefore $(\xi_P,\eta_P)\in\cU$. Hence there is an $(\al_j,\be_j)$, which we shall denote by $(\al_P,\be_P)$, such that $(\xi_P,\eta_P)\in\cO(\al_P)\times\cO(\be_P)\subset\cU$. \pa

$\al_P$ and $\be_P$, as preorders, are weakenings of $\xi_P$ and $\eta_P$, resp. Similarly to what we had in the proof of Theorem \ref{Thm:thm1}, also here for each link $\cI\!\cE$ in the chain of strict inequalities (which ``follow'' the staircase) in the definition of $\xi_P$ and $\eta_P$, flipping that inequality will create a full preorder $\xi'$ or $\eta'$ such that with it replacing $\xi_P$ or $\eta_P$, resp., one stair will be moved, and the intersection of the crosses will no longer be empty (i.e.\ the union of their complements will no longer be the whole $R_0$). Therefore $(\xi',\eta_P)$ or $(\xi_P,\eta')$, resp.\ is not in $\cU$, thus is not in $\cO(\al_P)\times\cO(\be_P)$, hence $\xi'$ (resp. $\eta'$) is not a strengthening of $\al_P$ (resp.\ $\be_P$), which compels $\al_P$ (resp.\ $\be_P$) to contain the inequality $\cI\!\cE$ with $\le$ instead of $<$. \pa

This means that for even $k$, $\xi_P$ contains the inequalities $b_{i_{k+1}}\le a_{i_k}$ for $1\le k\le N-1$ and $a_{i_k}\le b_{i_{k+3}}$ for $1\le k\le N-3$, while $\eta_P$ contains the inequalities $d_{i_k}\le c_{i_{k-1}}$ for $2\le k\le N$ and $c_{i_{k-3}}\le d_{i_{k}}$ for $4\le k\le N$. \pa

Now, again as in the proof of Theorem \ref{Thm:thm1}, one considers another permutation $P'=(i'_1,\ldots,i'_N)$ and defines $k\mapsto\ell=\ell(k)$ by $i_k=i'_{\ell(k)}$.
Thus {\em $\ell$ maps even numbers into even numbers and odd numbers into odd numbers}.
We shall prove that $\al_{P'}=\al_P$ is possible for at most four $P'$'s, i.e.\ $P\mapsto\al_P$ is at most four-to-one.
\pa

Indeed, suppose first that there exist $k,k'=1,\ldots,N$, with $k$, hence $\ell(k)$, even and $k'$, hence $\ell(k')$, odd, such that $k+3\le k'$ and $\ell(k')\le\ell(k)+1$. Then in $\xi_{P'}$,\,\, $a_{i_k}=a_{{i'}_{\ell(k)}}>b_{{i'}_{\ell(k)+1}}\ge b_{{i'}_{\ell(k')}}=b_{i_{k'}}$. Thus $\al_{P'}$, a weakening of $\xi(P')$, cannot contain $a_{i_k}\le b_{i_{k'}}$, which $\al_P$ contains, and one finds that in this case $\al_{P'}\neq\al_P$. Consequently, if $\al_{P'}=\al_P$ then
\BE\label{eq:kIIal}
k'\ge k+3\Rightarrow\ell(k')\ge\ell(k)+3,\quad\text{if }k,\ell(k)
\text{ even},k',\ell(k')\text{ odd}.
\EE
Then for any even $k$, the set of odds $\ge k+3$ must map by $\ell$ into the set of odds $\ge\ell(k)+3$. Consequently $\ell(k)\le k$ if $k\le N-4$. But since one can plug here $\ell^{-1}$ for $\ell$, also $k\le\ell(k)$ if $\ell(k)\le N-4$. This compels $\ell$ to fix all evens, except possibly interchanging the two evens among $N-3,N-2,N-1,N$. Also, for any odd $k'$, the set of evens $\le k'-3$ must map into the set of evens $\le\ell(k')-3$, hence $\ell(k')\ge k'$ if $k'\ge5$, and, plugging $\ell^{-1}$ for $\ell$, also $k'\ge\ell(k')$ if $\ell(k')\ge5$. This compels $\ell$ to fix all odds, except possibly interchanging $1$ and $3$. \pa

We have found that $P\mapsto\al_P$ is at most four-to-one, and there are at least $\LP\lfloor\half N\rfloor!\RP^2$ permutations $P$. Hence $M\ge\tfr14\LP\lfloor\half N\rfloor!\RP^2$ and the number of queries must be
$$m\ge2\log_6\LP\half\lfloor\half N\rfloor!\RP=\Om(N\log N).$$
\end{Prf}
\begin{Rmk}\label{Rmk:rmk2}
As in Setting I, but here less trivially (see Remark \ref{Rmk:rmk1}), the configuration that we used in the proof of Theorem \ref{Thm:thm2} (Figures \ref{fig:fig1} and \ref{fig:fig2}) provides an example where $N$ ``crosses'' have empty intersection while any proper subset of them does not. In the language of piercing pairs of intervals on two lines, we have an example of $N$ pairs of intervals, for any proper subset of them there is a piercing pair of points, while for the whole set there is none. Thus there is no ``Helly property'' here.
\end{Rmk}

\section{Outline of an $O(N\log N)$-steps algorithm for the piercing problem}

We consider setting II in dimension $2$ (the piercing problem) and speak in the language of ``crosses''. We assume we already know the ordering of $X$ and $Y$ (these are the sets of all the endpoints of the sides of the crosses), which, as is well known, can be found in $O(N\log N)$ steps (see the appendix). Hence we may assume that these endpoints are some integers between $0$ and $2N$, with the geometric ordering coinciding with the natural ordering of integers.
\footnote{We thank Amir Yehudayoff for instructive correspondence in the way to find such as algorithm.} \pa

Our aim is the intersection of all crosses. Now, every cross is the complement of four ``corner'' rectangles, which we refer to as the north-west, north-east etc. The intersection of the crosses is the complement of the union of all $4N$ corner rectangles, and we shall separate this union to four parts: the union of the north-west corners, the union of the north-east corners, etc. We shall construct each of the these four unions and then check the complement of all of them. \pa

For example, construct the union of the north-west rectangles: it is the area above an increasing ``function'' consisting of horizontal and vertical intervals. We construct this ``function'' in steps, storing it by pointers at its vertices (which we can think of as pairs of integers), pointing to the next vertex. The ``function'' will be successively amended by considering each next north-west ``corner'' (i.e.\ that of the next cross). For each such amendment -- adding one more corner -- we have to find out where the two arms of the new corner intersect the already existing increasing ``function'' (this requires $O(\log N)$ steps by the halving method, carried out on the set $\{0,1 ,\ldots,2N\}$) and then amend ($O(1)$ steps). Thus the union of the north-west corners is constructed in $O(N \log N)$ steps and consists of at most $2N$ (horizontal or vertical) intervals (at most two from each corner). \pa

Now we have four ``monotone functions'', so the intersection of the crosses is the area below the minimum of an increasing and a decreasing such ``function'' and above
the maximum of another increasing and another decreasing ``function''. The minimum
and maximum are constructed in $O(\log N)$ steps (needed to find the intersection of the two functions by halving) and the intersection of the crosses is non-empty if and only if the ``upper'' function is not below the ``lower'' one, and if we find a place where it is not below, then we have found a point in the intersection of the crosses. To check these
we follow the two ``functions'' from left to right, which takes $O(N)$ steps. \pa

Note that this leaves open the similar issue for setting I in $2$ dimensions - to check, in  $O(N\log N)$ steps, whether the union of $N$ rectangles covers the whole rectangle. \pa

\appendix\section{Ordering $N$ elements by $O(N\log N)$ queries}

The well-known procedure is as follows.\pa

Assume $N=2^n$ and suppose we have already enough information to order two given halves
$A$ and $B$ of the set. In order to finish ordering the whole set take the smallest element $a_1$ in $A$ in $A$'s now known ordering. Query to compare it with the smallest element of $B$ in its known ordering, then with the second, etc.\ until one finds the smallest $b_k\in B$ so that $a_1<b_k$. Then take the next element $a_2$ of $A$ and compare it with the elements of $B$ from $b_k$ onward. In $N$ queries we will check the whole $A$ and $B$ and thus know the ordering of the whole set. Similarly, to order $A$ and $B$ we need orderings of four quarters, thus $\frac12N+\frac12N=N$ queries. To order the quarters we need orderings of eight one-eights and $N$ queries. Hence by $nN=O(N\log N)$ queries the whole set can be ordered.


\end{document}